\documentclass[10pt,a4paper]{article}
\usepackage[utf8]{inputenc}
\usepackage{amsfonts}
\usepackage{amssymb}
\usepackage{amsmath}
\usepackage{color}
\usepackage{graphicx}
\usepackage{placeins}
\usepackage{amsthm}
\usepackage{enumerate}
\PassOptionsToPackage{normalem}{ulem}
\usepackage{ulem}
\usepackage[unicode=true,colorlinks,citecolor=linkcolor,linkcolor=linkcolor,urlcolor=urlcolor]{hyperref}
\usepackage[left=3cm,right=3cm,top=2.5cm,bottom=2.5cm]{geometry}
\usepackage{array}
\usepackage{cases}

\usepackage{color}
\usepackage{dsfont}

\theoremstyle{plain}
\newtheorem{theorem}{\protect Theorem}[section]
\newtheorem{proposition}[theorem]{\protect Proposition}
\newtheorem{definition}[theorem]{\protect Definition}
\newtheorem{remark}[theorem]{\protect Remark}
\newtheorem{assumption}[theorem]{\protect Assumption}

\definecolor{linkcolor}{rgb}{0,0,0.502}
\definecolor{urlcolor}{rgb}{1,0,0}

\begin{document}

\title{Stimulus sensitivity of a spiking neural network model}

\author{Julien Chevallier\footnote{{\sf{e-mail: \href{mailto:julien.chevallier1@univ-grenoble-alpes.fr}{julien.chevallier1@univ-grenoble-alpes.fr}}}}\\
AGM, UMR-CNRS 808, Université de Cergy-Pontoise, France \\
Univ. Grenoble Alpes, CNRS, LJK, 38000 Grenoble, France
}

\date{}

\maketitle

\begin{abstract}
Some recent papers relate the criticality of complex systems to their maximal capacity of information processing. In the present paper, we consider high dimensional point processes, known as age-dependent Hawkes processes, which have been used to model spiking neural networks. Using mean-field approximation, the response of the network to a stimulus is computed and we provide a notion of \emph{stimulus sensitivity}. It appears that the maximal sensitivity is achieved in the sub-critical regime, yet almost critical for a range of biologically relevant parameters.
\end{abstract}
 \smallskip\\
\textit{Physics and Astronomy Classification}: 87.18.Sn \and 87.19.La \and 05.70.Jk

\section{Introduction}

Self organization near critical points has been highlighted for a wide range of complex systems such as forest fires \cite{malamud1998forest} or computer networks \cite{valverde2002self} (see \cite{bak1991self,jensen1998self} for reviews). A great number of recent papers relate criticality of complex systems with optimal transmission of information for gene expression dynamics \cite{nykter2008gene}, swarm dynamics \cite{vanni2011criticality} and especially neural networks \cite{beggs2008criticality,kinouchi2006optimal,larremore2011predicting,shriki2016optimal}.

In \cite{kinouchi2006optimal}, Kinouchi and Copelli suggest to quantify the transmission of information in their neural network model thanks to the \emph{dynamic range} (expressed in decibel and mainly used for sound signals). It measures the range of stimuli for which a perturbation of the stimulus can be accurately described via a perturbation of the system's response. They consider a discrete time model where the stimulus is represented by a Bernoulli process (that is the time discrete version of the Poisson process) and the interactions between neurons follow an Erd\H{o}s-R\'enyi graph. For such a model, the parameter driving the criticality is the \emph{average branching ratio} (the branching ratio of each neuron measures its excitatory strength). A mean-field approximation of this discrete time stochastic model, supported by numerical simulations, shows that the dynamic range is first increasing then decreasing as a function of the average branching ratio and the maximum is achieved at criticality.

The neural network model used in the present article is known as age dependent Hawkes processes \cite{chevallier2017mean}: the spike train of each neuron, that is the set of the times at which it fires, is described by a point process. It can be thought as the time continuous version of the Kinouchi-Copelli model. The stimulus is described as a Poisson process and the interactions are described by a delay kernel function. Instead of the dynamic range, we introduce the notion of \emph{stimulus sensitivity} to quantify information processing. This notion is local compared to the dynamic range: it depends on the stimulus strength and the spontaneous activity of the neurons. The mean-field approximation of age dependent Hawkes processes yields a system of partial differential equations (PDE) as proved in \cite{chevallier2017mean}. Via analytical computations made on the stationary version of this PDE system, we show that the sensitivity is first increasing then decreasing as a function of our relevant parameter which quantifies the average connectivity strength. The main result states that the parameter value achieving the maximum is strictly less than the critical value (Theorem \ref{thm:1}).

It is worth to mention here the recent paper \cite{cassandro2017information} where the question of optimal transmission of information is addressed via purely analytical results (no numerical simulations) and without any mean-field approximation. The model considered there is the contact process and the sensitivity is measured with respect to an initial stimulus (instead of a constant input like in \cite{kinouchi2006optimal} or in the present paper) as a total variation distance. Using duality and coupling they are able to prove that maximal sensitivity is achieved in the super-critical regime in contrast to our main result. Optimal transmission seems then to be achieved below or above the critical point depending on the model. Nevertheless, optimum is linked with near-critical systems for biologically relevant parameter ranges (see Theorem \ref{thm:1} and the paragraph thereafter).

The paper is organized as follows. In Section \ref{sec:notation:definition}, the model is introduced and the main result is stated. Then, the response of the system to the stimulus is computed on the stationary version of the PDE as an explicit function of the model parameters in Section \ref{sec:steady:state:computations} and the main result is proved in Section \ref{sec:proof}.

\section{Model definition and main result}
\label{sec:notation:definition}

We consider here a system composed of age-dependent Hawkes processes with mean-field interaction (see \cite{chevallier2017mean} for more insight on this model). The parameters of the model are:
\begin{itemize}
\item a positive integer $n$ which is the number of neurons in the network;
\item an intensity function $\Psi:\mathbb{R}_{+}\times \mathbb{R} \rightarrow \mathbb{R}_{+}$;
\item an interaction function $h:\mathbb{R}_{+}\to \mathbb{R}$;
\item and a distribution $Q$ underlying the matrix of independent and identically distributed synaptic weights $(\alpha_{ij})_{i,j=1,\dots,n}$.
\end{itemize}

For each $i$ in $\{1,\dots,n\}$, the spike train of neuron $i$ is described by $N^{i}$ which is a point process with intensity $(\lambda^{i}_{t})_{t\geq 0}$ given by $\lambda^{i}_{t}= \Psi(S^{i}_{t-}, X^{i}_{t})$, where $S^{i}_{t-}$ is the \emph{age variable} and $X^{i}_{t}$ is the \emph{interaction variable} respectively defined by
\begin{equation*}
S^{i}_{t-}:=t-\sup\{T\in N^{i}, T<t\} \ \text{ and } \  X^{i}_{t}:= \frac{1}{n} \sum_{j=1}^{n} \int_{0}^{t-} \alpha_{ij} h(t-t') N^{j}(dt').
\end{equation*}
The age variable counts the time elapsed since the last spike of the neuron whereas the interaction variable sums up the influence of all spikes onto a given neuron at time $t$: the function $h$ deals with the delay in the transmission of information and $\alpha_{ij}$ represents the strength of the influence of spikes of neuron $j$ onto neuron $i$.

\begin{assumption}\label{ass:main:frame}
There exist $\mu,\delta\geq 0$, such that $\Psi(s,x)= (\mu+x) \mathbf{1}_{s\geq \delta}$, that is $\Psi(s,x)= \mu+x$ if $s\geq \delta$ and $\Psi(s,x)= 0$ otherwise. For all $t\geq 0$, $h(t)\geq 0$, $\int_{0}^{t} h(t')^{2} dt'<+\infty$ and $\int_{0}^{+\infty} h(t') dt'=1$. The distribution $Q$ is integrable, namely $\mathbb{E}\left[ |\alpha_{ij}| \right]<+\infty$.
\end{assumption}

Under Assumption \ref{ass:main:frame}, $\mu$ can be interpreted as the constant input signal received by the system, or as the spontaneous activity of the neurons, or as the sum of the two. The parameter $\delta$ is the duration of the strict refractory period: the delay between two spikes of any neuron is larger than $\delta$. In particular, the mean activity of one neuron cannot exceed $1/\delta$ even if the spontaneous rate $\mu$ is arbitrarily large (see for instance Figure \ref{fig:activity}). When the rate is close to the upper-bound $1/\delta$, we say that the system is close to saturation.

The parameter which controls the qualitative long time behaviour here is $\alpha:= \mathbb{E}\left[ \alpha_{ij} \right]$. It represents the average connectivity strength: for instance, if $Q$ is a Bernoulli distribution then the interaction graph is of Erd\H{o}s-R\'enyi type and $\alpha$ is the probability of connection between two neurons. In particular, the value of $\alpha$ tells if the system is sub-critical or super-critical. The critical value is $\alpha_{c}:=1$. Indeed, in the $\delta=0$ case, the system goes to a finite steady state activity when $\alpha<\alpha_{c}$ whereas the activity grows to infinity when $\alpha\geq \alpha_{c}$.

It is proved in \cite[Corollary 4.5.]{chevallier2017mean} that such a mean field complex system is well approximated, when $n$ goes to infinity, by the following non linear partial differential equation with parameters $(\mu,\alpha,\delta)$,
\begin{equation}\label{eq:PPS}
\begin{cases}
\displaystyle \left(\frac{\partial}{\partial t}+\frac{\partial}{\partial s}\right) u(t,s) + \Psi(s,x(t)) u(t,s) =0,\\
\displaystyle u(t,0)= \int_{0}^{+\infty} \Psi(s,x(t)) u(t,s) ds,\vspace{0.5em}\\
x(t) = \alpha \int_{0}^{t}  h(t-t') u(t',0) dt'.
\end{cases}
\end{equation}
In the system above, $u(t,s)$ represents the probability density of the age $s$ at time $t$ of any representative neuron in the asymptotic $n\to +\infty$. The variable $x(t)$ is the deterministic counterpart of the stochastic interaction variables $X^{i}_{t}$. As highlighted by the second equation of the PDE system \eqref{eq:PPS}, the network average instantaneous activity is given by $a_{t}:=u(t,0)$.

The purpose of the present paper is then to consider the stationary regime\footnote{In fact, exponential convergence of the solution of the PDE system \eqref{eq:PPS} to the unique solution of its stationary version \eqref{eq:stationary:PPS} can be proved in the weak ($\alpha<\alpha_{\rm weak}$) and the strong ($\alpha>\alpha_{\rm strong}$) connectivity regime \cite{weng2015general}. However, the result are not quantitative: $\alpha_{\rm weak}$ and $\alpha_{\rm strong}$ are unknown. We do not know how they compare to $\alpha_{c}$.} characterized by
\begin{equation}\label{eq:stationary:PPS}
\begin{cases}
\displaystyle \frac{\partial}{\partial s} u_{\infty}(s) + \Psi(s,x_{\infty}) u_{\infty}(s) =0,\\
\displaystyle a_{\infty}:=u_{\infty}(0)= \int_{0}^{+\infty} \Psi(s,x_{\infty}) u_{\infty}(s) ds,\\
\displaystyle x_{\infty} = \alpha \,a_{\infty}.
\end{cases}
\end{equation}
Under Assumption \ref{ass:main:frame}, the stationary solution $u_{\infty}$ is unique - see Section \ref{sec:steady:state:computations}. Furthermore, the steady state activity $a_{\infty}$ can be expressed as a function of the parameters $(\mu,\alpha,\delta)$ as stated below in Proposition \ref{prop:a:infty:expression} and the sensitivity is defined as follows.

\begin{definition}
The \emph{stimulus sensitivity} of the system with parameters $(\mu,\alpha,\delta)$ is denoted by $\sigma$ and defined as
\begin{equation*}
\sigma=\sigma(\mu,\alpha,\delta) := \frac{\partial}{\partial \mu} a_{\infty}(\mu,\alpha,\delta).
\end{equation*}
\end{definition}
This quantity measures the sensitivity of the system to a perturbation of a given input signal $\mu$. It can also be interpreted as the sensitivity to the occurrence of a small input signal onto a system with spontaneous activity equal to $\mu$. The main result of the paper then reads as follows.

\begin{theorem}\label{thm:1}
For every given $\mu$, $\delta>0$, there exists a unique $\alpha_{m}:=\alpha_{m}(\mu,\delta)\geq 0$ for which the sensitivity $\sigma$ is maximal. Furthermore, $\alpha_{m}$ is non-increasing with respect to the product $\mu\delta$ between the two extremal behaviours:
\begin{equation}\label{eq:results:thm:1}
\begin{cases}
\text{if } \mu\delta\geq 1/2 \text{ then } \alpha_{m}=0,\\
\text{if } \mu\delta\to0 \text{ then } \alpha_{m}\to\alpha_{c}:=1.
\end{cases}
\end{equation}
In particular, $\alpha_{m}<\alpha_{c}$, which means that maximal sensitivity is reached below criticality.
\end{theorem}

\begin{remark}
The results stated above depend on the dimensionless\footnote{It is dimensionless since $\mu$ is a frequency and $\delta$ a duration.} product $\mu\delta$. This is a consequence of the fact that the dependence of $\sigma$ with respect to (w.r.t.) $(\mu,\delta)$ is reduced to a dependence w.r.t. $\mu\delta$ as it appears in Proposition \ref{prop:a:infty:limits} below. 
\end{remark}

\begin{figure}[h]
\begin{center}
\includegraphics[width=0.55\textwidth]{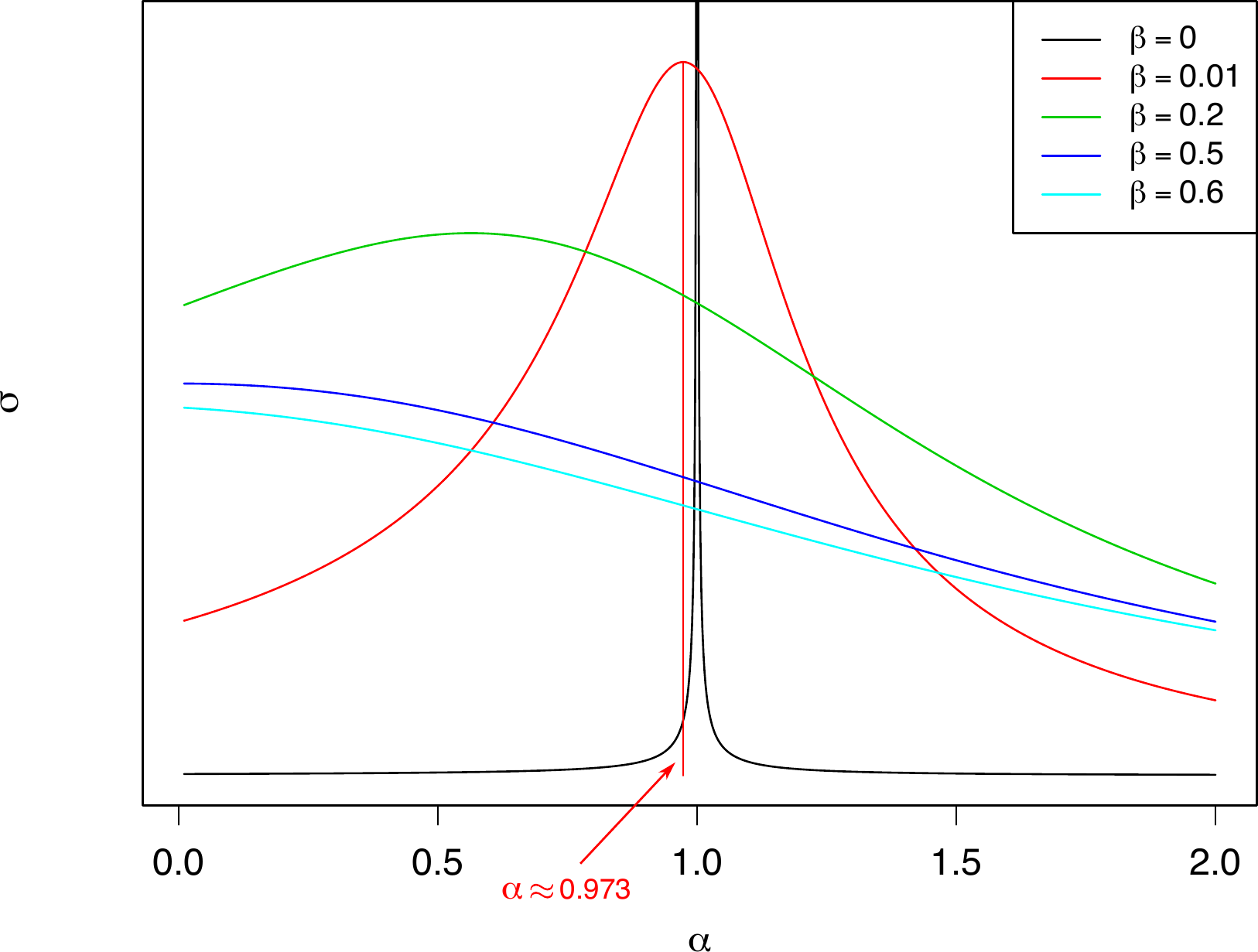}
\end{center}
\caption{\label{fig:sensitivity}
Sensitivity $\sigma$ as a function of $\alpha$ for different values of $\beta:=\mu\delta\in\{0,0.01,0.2,0.5,0.6\}$. The vertical scale differs from one curve to the other : only the shapes of the curves are informative. The red vertical line corresponds to the optimal sensitivity when $\beta=0.01$ so $\alpha_{m}(\mu,\delta)\approx 0.973$ when $\mu\delta=0.01$.
}
\end{figure}

Let us give some intuition on the extremal behaviours \eqref{eq:results:thm:1}. When $\mu\delta$ is large, the system is close to saturation and so it is little sensitive to perturbations. In that case, the unconnected regime ($\alpha=0$) is further from saturation and the system has more latitude to sense perturbations. The case $\mu\delta\to 0$ corresponds to a system of neurons with low spontaneous activity and/or driven by very weak inputs. Then, the system is far from saturation in the unconnected regime and optimal sensitivity is reached near criticality. This is consistent  This second case corresponds to biologically relevant parameters range: the spontaneous activity is around one Hertz whereas the refractory period is around one millisecond. For instance, with $\mu=2$ Hz and $\delta=5$ ms, we have $\mu\delta=0.01$ and numerical computations give $\alpha_{m}>0.97$ (see Figure \ref{fig:sensitivity}).

\section{Steady state activity}
\label{sec:steady:state:computations}

The following proposition gives the expression of the steady activity $a_{\infty}$ as a function of the parameters.

\begin{proposition}\label{prop:a:infty:expression}
For all $\mu>0$,
\begin{equation}\label{eq:a:infty:expression}
a_{\infty}(\mu,\alpha,\delta)=
\left\{
\begin{aligned}
& \mu/(1-\alpha) & & \text{if $\delta=0$ and $\alpha<\alpha_{c}$},\\
& 1/(\delta+1/\mu) & & \text{if $\delta\geq 0$ and $\alpha=0$},\\
& \frac{1}{\delta} - \frac{1+\alpha+\mu\delta - \sqrt{\Delta(\mu,\alpha,\delta)}}{2\alpha\delta} & & \text{if $\delta>0$ and $\alpha>0$},\\
\end{aligned}
\right.
\end{equation}
where $\Delta(\mu,\alpha,\delta):= (1+\mu\delta-\alpha)^{2} + 4\mu\alpha\delta$.
\end{proposition}
Let us elaborate on the different cases involved in the expression \eqref{eq:a:infty:expression} of the steady activity $a_{\infty}$. The first line corresponds to standard Hawkes processes: it is well known that the activity raises up to a finite value when $\alpha<\alpha_{c}$ whereas it goes to infinity when $\alpha\geq \alpha_{c}$\footnote{This can be proved thanks to the cluster representation of linear Hawkes processes \cite{hawkes_1971} and similar results available on Galton-Watson trees.}. The second line corresponds to the framework of an unconnected network. In that case, the underlying stochastic processes are independent Poisson processes with dead time \cite{Reimer2012}. Looking closely to \eqref{eq:a:infty:expression}, we retrieve the already mentioned fact that the activity is upper-bounded by $1/\delta$: they can only saturate up to their refractoriness. Notice that they achieve saturation in the limit $\mu\to +\infty$ as highlighted below.

\begin{figure}[h]
\begin{center}
\includegraphics[width=0.49\textwidth]{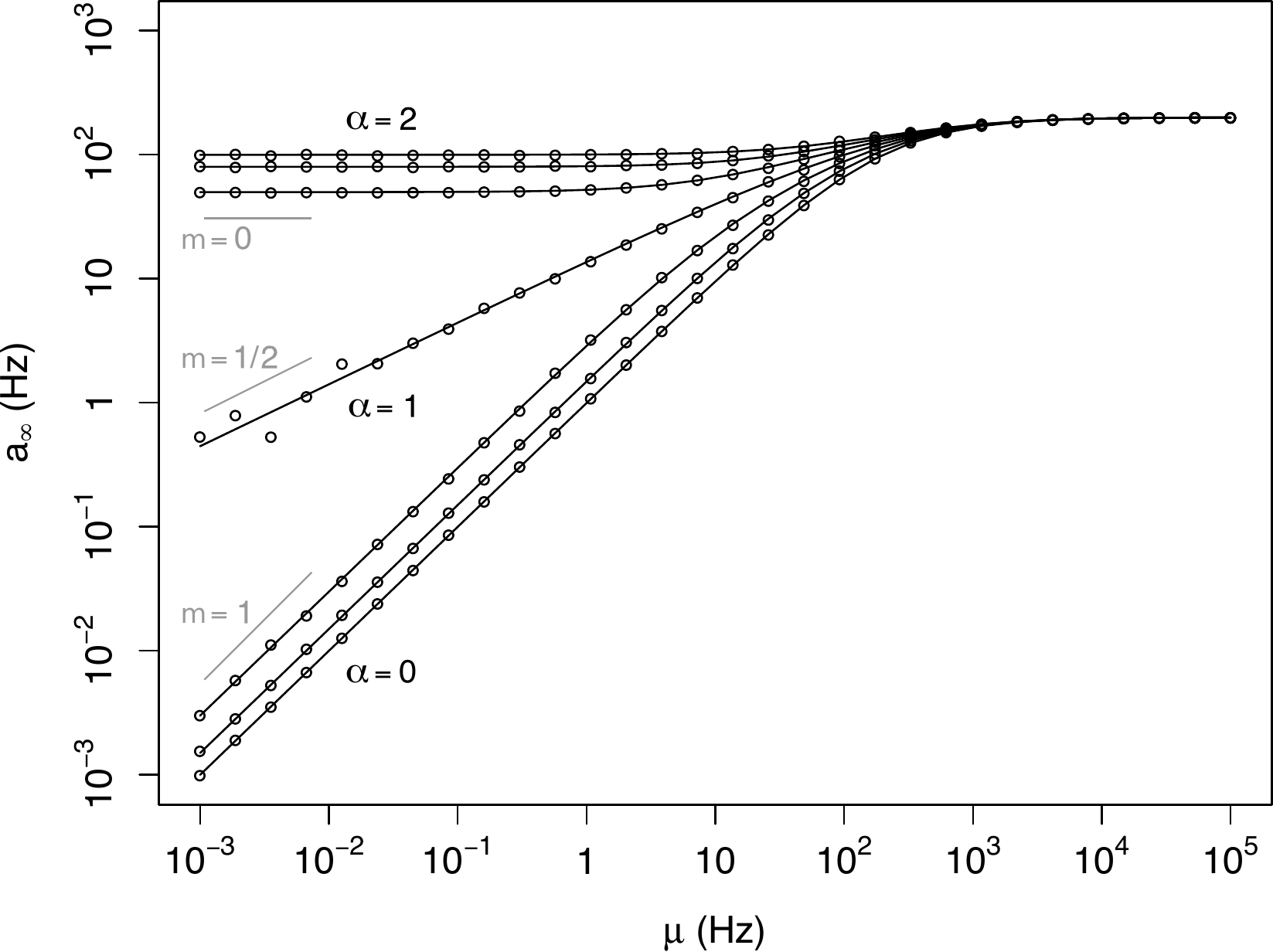}
\includegraphics[width=0.49\textwidth]{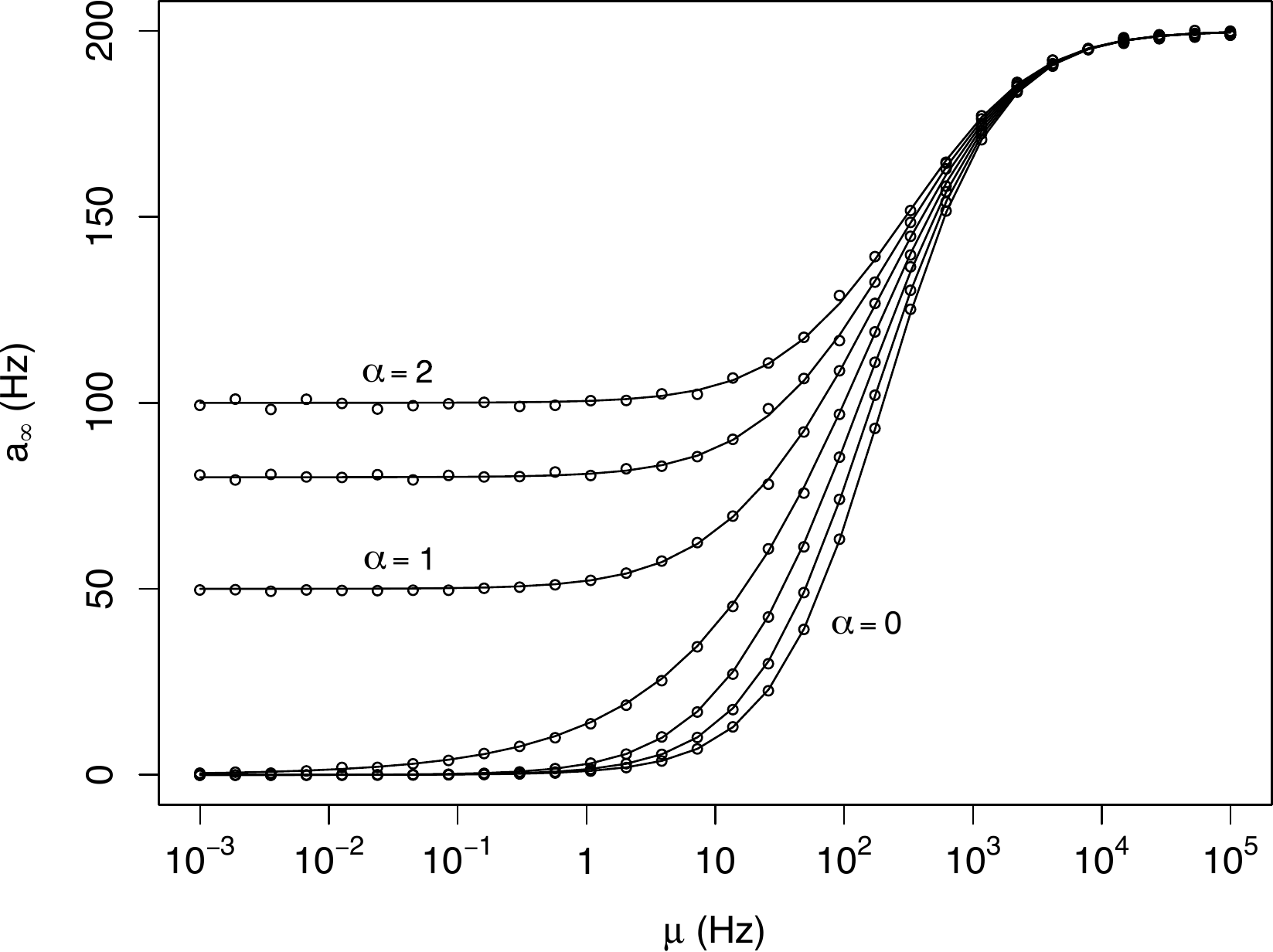}
\end{center}
\caption{\label{fig:activity}
Steady activity $a_{\infty}$ as a function of $\mu$ for different values of $\alpha$ between $0$ and $2$ (exact values are $k/3$, $k=0,\dots,6$) and fixed $\delta=5$ ms. The curves are obtained via the formula stated in Proposition \ref{prop:a:infty:expression} and the circles are obtained through numerical simulations of age-dependent Hawkes processes with $n=1000$ neurons and the activity is averaged over the $5000$ first spikes produced by the system. Left: both axes are in log scale. For $\alpha<1$, the slope is $m=1$ since the activity is linear w.r.t. $\mu$. For $\alpha=1$, $m=1/2$ since the activity is of square-root order - see the Taylor expansion \eqref{eq:taylor:expansion}. For $\alpha>1$, the slope is null since the activity without input is non-zero - see \eqref{eq:a:infty:limits}.
Right: the vertical axis is in linear scale.
}
\end{figure}

The following proposition gives an useful expression for the sensitivity and states the following expected result: in response to a larger input, the system shows a larger activity.
\begin{proposition}\label{prop:a:infty:limits}
For all $\mu,\alpha,\delta>0$, 
\begin{equation}\label{eq:sigma:expression}
\sigma(\mu,\alpha,\delta) = -\frac{1}{2\alpha} + \frac{1+\mu\delta+\alpha}{2\alpha\sqrt{\Delta(\mu,\alpha,\delta)}},
\end{equation}
where $\Delta$ is defined in Proposition \ref{prop:a:infty:expression}. In particular, 
\begin{equation}\label{eq:limit:sigma:mu:to:0}
\lim_{\mu\delta\to 0} \sigma(\mu,\alpha,\delta)=
\left\{
\begin{aligned}
& \frac{1}{1-\alpha} & & \text{if $\alpha<1$},\\
& \frac{1}{\alpha(\alpha-1)} & & \text{if $\alpha>1$},\\
\end{aligned}
\right.
\end{equation}
and $a_{\infty}$ is increasing w.r.t. $\mu$ and remains between the two following bounds,
\begin{equation}\label{eq:a:infty:limits}
\begin{cases}
\lim_{\mu\to 0} a_{\infty}(\mu,\alpha,\delta) = \max\big( (\alpha-1)/(\alpha\delta) ,0 \big),\\
\lim_{\mu\to +\infty} a_{\infty}(\mu,\alpha,\delta)= 1/\delta.
\end{cases}
\end{equation}
\end{proposition}
Looking closely to the expression \eqref{eq:sigma:expression} of the sensitivity, we notice that it does not really depend on the full couple $(\mu,\delta)$ but only on the dimensionless product $\mu\delta$ which is therefore non sensitive to a time change (which could be expected since we consider steady states).
Due to the limit expression \eqref{eq:limit:sigma:mu:to:0}, we expect that optimal sensitivity is achieved near criticality in the low spontaneous activity regime ($\mu\delta\to 0$). Furthermore, the sensitivity goes to infinity in that limit case and the steady activity has a square-root behaviour (instead of a linear one) w.r.t. $\mu$ as stressed by the following Taylor expansion,
\begin{equation}\label{eq:taylor:expansion}
a_{\infty}(\mu,\alpha,\delta)= \frac{1}{\sqrt{\delta}} \mu^{1/2} + o(\mu^{1/2}) \quad \text{ when $\alpha=\alpha_{c}$.}
\end{equation}
The first line of \eqref{eq:a:infty:limits} shows a transition at $\alpha=\alpha_{c}$. When $\alpha>\alpha_{c}$, the response to an arbitrary small input is not arbitrary small - see the right plot of Figure \ref{fig:activity}.

\begin{proof}[Propositions \ref{prop:a:infty:expression} and \ref{prop:a:infty:limits}]
From the first equation of the PDE system \eqref{eq:stationary:PPS}, we deduce that $u_{\infty}(s) = a_{\infty} e^{-(\mu+x_{\infty})\max(s-\delta,0)}$.
Since $\mu>0$, there is still mass conservation for large times\footnote{For all $t\geq 0$, $u(t,\cdot)$ remains a probability density and so it is for $u_{\infty}$.} so
\begin{equation*}
1= \int_{0}^{+\infty} u_{\infty}(s)ds = a_{\infty} \left(\delta+ \frac{1}{\mu+x_{\infty}}\right).
\end{equation*}
From this, we easily deduce \eqref{eq:a:infty:expression}. 

Then, the expression of $\sigma$ is deduced from the expression \eqref{eq:a:infty:expression} of $a_{\infty}$ by differentiating w.r.t. $\mu$. Since we can rewrite $\Delta(\mu,\alpha,\delta)=(1+\mu\delta+\alpha)^{2} - 4\alpha$, we notice that $\sigma$ is positive and so $a_{\infty}$ is increasing. Finally, the limits \eqref{eq:a:infty:limits} also follow from the last line of \eqref{eq:a:infty:expression}.
\end{proof}

\section{Proof of Theorem \ref{thm:1}}
\label{sec:proof}

We start with the existence and uniqueness of $\alpha_{m}$ by studying the derivative of $\sigma$ namely
\begin{equation}\label{eq:d:alpha:sigma:expression}
g(\mu,\alpha,\delta):=\frac{\partial}{\partial\alpha}\sigma(\mu,\alpha,\delta) = \frac{\Delta^{3/2} - (1+\mu\delta)\Delta - \alpha\big((\mu\delta+\alpha)^{2}-1\big)}{2\alpha^{2}\Delta^{3/2}},
\end{equation}
where $\Delta=\Delta(\mu,\alpha,\delta)$ is given in Proposition \ref{prop:a:infty:expression}. Let us denote $\beta=\mu\delta$, $\Delta=\Delta(\alpha,\beta):=(1+\beta-\alpha)^{2}+4\alpha \beta$ and $g(\alpha,\beta)=g(\mu,\alpha,\delta)$ in the following. First, let us remark that
\begin{equation}\label{eq:lim:derivative:alpha:to:0}
\lim_{\alpha\to 0} g(\alpha,\beta) = \frac{1-2\beta}{(1+\beta)^{4}}.
\end{equation}
Furthermore, let $f(\alpha,\beta)$ denote the numerator of the derivative $g(\alpha,\beta)$ and $f^{c}(\alpha,\beta)$ denote its conjugate quantity namely
\begin{equation*}
f^{c}(\alpha,\beta):=\Delta^{3/2} + (1+\beta)\Delta +\alpha((\beta+\alpha)^{2}-1).
\end{equation*}
The product of $f$ and $f^{c}$ gives
\begin{equation*}
\Delta^{3} - \left[(1+\beta)\Delta +\alpha((\beta+\alpha)^{2}-1)\right]^{2}=-4\alpha^{2} P(\alpha,\beta),
\end{equation*}
with $P(\alpha,\beta):=2\alpha^{3}+(6\beta-5)\alpha^{2}+(6\beta^{2}-6\beta+4)\alpha+2\beta^{3}+3\beta^{2}-1$.

When $\beta\geq 1/2$, it is easy to check that $f^{c}(\alpha,\beta)\geq 0$ for all $\alpha>0$, so that the sign of $g(\alpha,\beta)$ is opposite to the sign of $P(\alpha,\beta)$. Yet $P(\beta,\alpha)\geq 2\alpha^{3}-2\alpha^{2}+(5/2) \alpha>0$ for every positive $\alpha$. Hence $\alpha_{m}=0$ as soon as $\beta\geq 1/2$ - first line of Equation \eqref{eq:results:thm:1}.

When $\beta<1/2$, we characterize $\alpha_{m}$ as the unique \emph{stationary point} of $\sigma$ (that is where its derivative is null) which is necessarily a zero of $P(\cdot,\beta)$. First, when $\alpha\geq \alpha_{c}=1$, we have $f^{c}(\alpha,\beta)\geq \Delta^{3/2}+(1+\beta)\Delta > 0$ so the sign of the derivative \eqref{eq:d:alpha:sigma:expression} is opposite to the sign of $P(\alpha,\beta)$. Yet $2\alpha^{3}-5\alpha^{2}+4\alpha-1\geq 0$, $2\beta^{3}+(3+6\alpha)\beta^{2}+6\alpha(\alpha-1)\beta>0$ and so $P(\beta,\alpha)>0$. Gathering this, we obtain that $f(\alpha,\beta)< 0$ as soon as $\alpha\geq 1$. Combined with limit value \eqref{eq:lim:derivative:alpha:to:0}, this ensures the existence of a stationary point for $\sigma$. Furthermore, if it is unique then it is necessarily a maximum.

Now, let us prove uniqueness by considering the zeros of the cubic polynomial $P(\cdot,\beta)$. Thanks to its discriminant, we obtain the three following cases depending on the value of $\beta$ w.r.t. $\beta_{0}:=1/27$ : 1/ if $\beta>\beta_{0}$ there is a unique real zero, 2/ if $\beta=\beta_{0}$ there are two distinct real zeros (multiplicities equal to one and two), 3/ if $\beta<\beta_{0}$ there are three distinct real zeros.

In the first case, it is clear that the stationary point of $\sigma$ is unique. In the other cases, we prove that the conjugate quantity $f^{c}$ has two zeros (up to multiplicity) between $\alpha=0$ and $\alpha=1$ thanks to the remark that $f^{c}(0,\beta)=2(1+\beta)^{3}>0$ and $f^{c}(1,\beta)> 0$. Indeed, $f^{c}(\alpha,\beta_{0})=0$ has a unique solution $\alpha=\alpha_{0}:=20/27$ with multiplicity two. Furthermore, it appears that
\begin{equation*}
f^{c}(\alpha_{0},\beta)=-A + B\beta + C\beta^{2} +\beta^{3} +(D+E\beta+\beta^2)^{3/2},
\end{equation*}
where $A,B,C,D,E$ are positive constants. In particular, it is increasing w.r.t. $\beta$ which implies that $f^{c}(\alpha_{0},\beta)< 0$ for all $\beta<\beta_{0}$ and so $f^{c}(\cdot,\beta)$ has two distinct zeros in that case. Hence the zero of $f$ and so the stationary point of $\sigma$ is unique and we have proved existence and uniqueness of $\alpha_{m}$.

Then, we compute
\begin{equation}\label{eq:derivative:f}
\frac{\partial}{\partial \beta} f(\alpha,\beta) = 3(1+\beta+\alpha)\Delta^{1/2} - 3(1+2\beta+\beta^{2}+2\alpha \beta+\alpha^{2}).
\end{equation}
Its conjugate quantity is clearly non negative so the sign of the quantity above is the same as the sign of
\begin{equation*}
9(1+\beta+\alpha)^{2}\Delta - 9(1+2\beta+\beta^{2}+2\alpha \beta+\alpha^{2})^{2} = -36\alpha^{2}.
\end{equation*}
Hence, $f(\alpha,\cdot)$ is decreasing for all $\alpha>0$. Yet, we know that if $\beta>1/2$ then $f(\alpha_{m}(\beta),\beta)=0$ and $\alpha_{m}(\beta)>0$ so that for all $\tilde{\beta}>\beta$, $f(\alpha_{m}(\beta),\tilde{\beta})<0$ and in particular $\alpha_{m}(\tilde{\beta})<\alpha_{m}(\beta)$ that is $\alpha_{m}$ is non-increasing.

Finally, let us denote $\overline{\alpha}:=\lim_{\beta\to 0} \alpha_{m}\leq 1$ which exists since it is non-increasing and bounded. Since $\sigma(\beta,\alpha_{m})\geq \sigma(\beta,1)\to +\infty$ when $\beta\to 0$ we have $\sigma(0,\overline{\alpha})=+\infty$. Yet, for all $\alpha<1$, $\sigma(\tilde{\alpha},\beta)\to 1/(1-\alpha)<+\infty$ when $\beta\to 0$ and $\tilde{\alpha}\to \alpha$. Hence $\overline{\alpha}=1$ which completes the proof.

\section{Conclusion and Outlooks}

Many studies relate optimization of functional properties (e.g. dynamic range or information transmission) of cortical networks with (near-)critical regime (see the review \cite{shew2013functional}). The present paper deals with stimulus sensitivity, which has to be related to dynamic range. Depending on the neural network model, optimal sensitivity seems to be achieved by critical \cite{kinouchi2006optimal}, or even super-critical \cite{cassandro2017information}, dynamics. In contrast, we have proved, analytically, that stimulus sensitivity of the mean field approximation of age-dependent Hawkes processes is maximized in the sub-critical regime yet near-critical for biologically relevant parameters values. This result shows nice agreement with estimation performed on \emph{in vivo} spiking activity (see the recent paper \cite{wilting2016branching} dealing with the issue of underestimation due to subsampling). Furthermore, lowering the input signal leads to optimization closer to criticality. This observation is consistent with the link between sensory deprivation and hallucinations \cite{shriki2016optimal}.

Hawkes processes and related models are known to exhibit event cascades in near critical regime (see for instance the revised Hawkes process in \cite{onaga2016emergence}). Up to our knowledge, results on the distribution of the cascade length are limited to upper-bounds on the distribution tail (exponential bounds are given in \cite{reynaud_2007}). It would be interesting to relate the near critical dynamics of Hawkes processes identified in this paper to power law distributed cascades observed in experiments \cite{arviv2015near,beggs2003neuronal}.

A possible extension to the present study is to consider time-varying inputs and especially periodic signals. Two questions then arise. 1) Can the system reproduce the periodic signal ? In other words, does the system \eqref{eq:PPS} exhibit periodic solutions under periodic stimulation ? 2) If so, what is the good notion of stimulus sensitivity and how does it relate to the critical point ? The major difficulty in that framework is that we need to handle two typical times : the period of the stimulation and the typical interaction time encoded in the interaction function $h$ (the latter was irrelevant in the present article since we only considered the stationary regime).

\paragraph{Acknowledgements}
The author would like to thank Eva Löcherbach for inspiring discussions on this topic. This research was supported by the project Labex MME-DII (ANR11-LBX-0023-01) and mainly conducted during the stay of the author at Universit\'e de Cergy-Pontoise.

\bibliographystyle{abbrv}
\bibliography{references}

\end{document}